\documentclass[preprint,12pt,1p]{elsarticle}
\makeatletter
\def\ps@pprintTitle{%
  \let\@oddhead\@empty
  \let\@evenhead\@empty
  \let\@oddfoot\@empty
  \let\@evenfoot\@oddfoot
}
\makeatother
\usepackage{graphics}
\usepackage{epsfig}
\usepackage{mathptmx}
\usepackage{amsmath}
\usepackage{amssymb}
\usepackage{hyperref}
\usepackage{fullpage}
\usepackage{amsthm}
\usepackage{geometry}
\usepackage{bigints}
\theoremstyle{plain}

\theoremstyle{definition}

\theoremstyle{example}

\theoremstyle{remark}

\newlength{\defbaselineskip}
\journal{Computer Aided Geometric Design}
\numberwithin{equation}{section}
\begin{document}
\begin{frontmatter}

\title{ Some definite integrals of Srinivasa Ramanujan and its consequences}
\author{M. I. Qureshi}
\author{*Showkat Ahmad Dar}
\ead{showkat34@gmail.com}
\address{miqureshi\_delhi@yahoo.co.in, showkat34@gmail.com}
\address{Department of Applied Sciences and Humanities , \\Faculty of  Engineering and Technology , \\Jamia Millia Islamia (A Central University), New Delhi, 110025, India.}
\cortext[cor1]{Corresponding author}
\begin{abstract}
In this paper, we obtain analytical solutions of some definite integrals of Srinivasa Ramanujan [Mess. Math., XLIV, 75-86, 1915] in terms of Meijer's $G$-function by using Laplace transforms of $ \sin(\beta x^{2}),\cos(\beta x^{2}), x\sin(\beta x^{2})$ and $x\cos(\beta x^{2})$. Further, we obtain some infinite summation formulas connected with Meijer's G-function and
numeric values of some infinite series. \\
\\
\textit{2010 AMS Classification:  44A10; 33E20; 33C60.}
 \end{abstract}

\end{frontmatter}
\section{Introduction and Preliminaries}
 For the sake of conciseness of this paper, we use the following notations \cite[p.33]{C3}\\
$~~~~~~~~~~~~~~~~~~~~~~~~~~~~\mathbb{N}:=\{1,2,...\};~~~~~~\mathbb{N}_{0}:=\mathbb{N}\cup\{0\};~~~~~~\mathbb{Z}_{0}^{-}:=\mathbb{Z}^{-}\cup\{0\}$,\\
and
\begin{equation}\label{cr}
\mathbb{R}_{>}:=\{a\in\mathbb{R} : a>0\} ~~~~; ~~~~\mathbb{C}_{>}:=\{p\in\mathbb{C} : \Re(p)>0\},
\end{equation}
where the symbols $\mathbb{N}$ and $\mathbb{Z}$ denote the set of natural number and integers; as usual, the symbols $\mathbb{R}$ and $\mathbb{C}$ denote the set of real and complex numbers.\\
The following three theorems of  Srinivasa Ramanujan are available in the collected papers \cite{H} edited by Hardy- Aiyar- Wilson, without giving any analytical proofs.
In our recent communication \cite{Q13} we have verified  numerically the following three theorems of Ramanujan, with the help of Wolfram Mathematica software[see \cite{Q13} in Tables 6.1, 6.2; 7.1,7.2;8.1,8.2].\\
\\
\\
\\
\textbf{Theorem I}: The first theorem of Ramanujan edited by  G.H. Hardy et.al \cite[p.59, eq.(3)and $(\acute{3})$]{H}and \cite[p.75, eq.(3 and $3^{'}$)]{R1} is given below\\
If
\begin{equation}\label{CR14}
\Phi_{1} (n)=\int_{0}^{\infty}\frac{\cos(\pi n x^{2})}{\cosh(\pi x)}dx,
\end{equation}
and
\begin{equation}\label{CR15}
\Psi_{1}(n)=\int_{0}^{\infty}\frac{\sin(\pi n x^{2})}{\cosh(\pi x)}dx,
\end{equation}
then
\begin{equation}\label{CR16}
\Phi_{1}(n)=\sqrt{\left(\frac{2}{n}\right)}~\Psi_{1}\left(\frac{1}{n}\right)+\Psi_{1}(n),
\end{equation}
and
\begin{equation}\label{CR17}
\Psi_{1}(n)=\sqrt{\left(\frac{2}{n}\right)}~\Phi_{1}\left(\frac{1}{n}\right)-\Phi_{1}(n),
\end{equation}
where $n\in\mathbb{R}_{>}$ .\\
\textbf{Theorem II}:The second theorem of Ramanujan edited by  G.H. Hardy et.al \cite[p.60, eq.(6) and $(\acute{6}) $]{H} and \cite[p.76, eq.(6 and $6^{'}$)]{R1} is given below\\
If
\begin{equation}\label{CR22}
\Phi_{2} (n)=\int_{0}^{\infty}\frac{\cos(\pi n x^{2})}{\{1+2\cosh\left(2\pi x/\sqrt{3}\right)\}}dx,
\end{equation}
and
\begin{equation}\label{CR23}
\Psi_{2}(n)=\int_{0}^{\infty}\frac{\sin(\pi n x^{2})}{\{1+2\cosh\left(2\pi x/\sqrt{3}\right)\}}dx,
\end{equation}
then
\begin{equation}\label{CR24}
\Phi_{2}(n)=\sqrt{\left(\frac{2}{n}\right)}~\Psi_{2}\left(\frac{1}{n}\right)+\Psi_{2}(n),
\end{equation}
and
\begin{equation}\label{CR25}
\Psi_{2}(n)=\sqrt{\left(\frac{2}{n}\right)}~\Phi_{2}\left(\frac{1}{n}\right)-\Phi_{2}(n),
\end{equation}
where $n\in\mathbb{R}_{>}$ .\\
\\
\textbf{Theorem III}:The  third theorem of Ramanujan edited by  G.H. Hardy et.al \cite[p.60, eq.(10)and $(10^{'})$;  p.66, after eq.(44)]{H}; see also \cite{B13} and \cite[p.76-77, eq.(10 and $10^{'}$)]{R1} is given below \\
If
\begin{equation}\label{CR26}
\Phi_{3} (n)=\int_{0}^{\infty}\frac{\cos(\pi nx)}{\{-1+\exp{(2\pi\sqrt{x})}\}}dx,
\end{equation}
and
\begin{equation}\label{CR27}
\Psi_{3}(n)=\frac{1}{2\pi n}+\int_{0}^{\infty}\frac{\sin(\pi nx)}{\{-1+\exp{(2\pi\sqrt{x})}\}}dx=\frac{1}{2\pi n}+\Psi_{3}^{*}(n),
\end{equation}
then
\begin{equation}\label{CR28}
\Phi_{3}(n)=\frac{1}{n}\sqrt{\left(\frac{2}{n}\right)}~\Psi_{3}\left(\frac{1}{n}\right)-\Psi_{3}(n),
\end{equation}
and
\begin{equation}\label{CR29}
\Psi_{3}(n)=\frac{1}{n}\sqrt{\left(\frac{2}{n}\right)}~\Phi_{3}\left(\frac{1}{n}\right)+\Phi_{3}(n)=\frac{1}{2\pi n}+\Psi_{3}^{*}(n),
\end{equation}
where $n\in\mathbb{R}_{>}$ .\\
For a particular value of $n=1$, in the Ramanujan's integrals (\ref{CR14}),(\ref{CR15}),(\ref{CR22}) and (\ref{CR23}),  we have
\begin{equation}\label{N1}
\Phi_1(1)=\int_{0}^{\infty}\frac{\cos(\pi x^{2})}{\cosh(\pi x)}dx=\frac{1}{2\sqrt{2}},
\end{equation}
\begin{equation}\label{N2}
\Psi_1(1)=\int_{0}^{\infty}\frac{\sin(\pi  x^{2})}{\cosh(\pi x)}dx=\frac{-1+\sqrt{2}}{2\sqrt{2}},
\end{equation}
\begin{equation}\label{N3}
\Phi_2(1)=\int_{0}^{\infty}\frac{\cos(\pi x^{2})}{\{1+2\cosh\left(2\pi x/\sqrt{3}\right)\}}dx=\frac{2-\sqrt{6}+\sqrt{2}}{8},
\end{equation}
\begin{equation}\label{N4}
\Psi_2(1)=\int_{0}^{\infty}\frac{\sin(\pi x^{2})}{\{1+2\cosh\left(2\pi x/\sqrt{3}\right)\}}dx=\frac{-\sqrt{12}+\sqrt{2}+\sqrt{6}}{8}.
\end{equation}
When we put $t=0$ in the results \cite[p.79, eq.23; p.80, eq.24 and p.84, eq.41]{R1}, after simplification we get above numeric values of the integrals (\ref{N1})-(\ref{N4}).\\
Again for particular values of $n$, some values of Ramanujan's integral $\Phi_{3} (n)$ \cite[p.85 (eq. 48)]{R1} are available.
\begin{equation}\label{RCG10}
\Phi_{3}(1)=\int_{0}^{\infty}\frac{\cos(\pi x)}{\{-1+\exp{(2\pi\sqrt{x})}\}}dx=\frac{2-\sqrt{2}}{8},
\end{equation}
\begin{equation}\label{RCG11}
\Phi_3(2)=\int_{0}^{\infty}\frac{\cos(2\pi x)}{\{-1+\exp{(2\pi\sqrt{x})}\}}dx=\frac{1}{16},
\end{equation}
\begin{equation}\label{RCG12}
\Phi_3(4)=\int_{0}^{\infty}\frac{\cos(4\pi x)}{\{-1+\exp{(2\pi\sqrt{x})}\}}dx=\frac{3-\sqrt{2}}{32},
\end{equation}
\begin{equation}\label{RCG13}
\Phi_3(6)=\int_{0}^{\infty}\frac{\cos(6\pi x)}{\{-1+\exp{(2\pi\sqrt{x})}\}}dx=\frac{13-4\sqrt{3}}{144},
\end{equation}
\begin{equation}\label{RCG14}
\Phi_3\left(\frac{1}{2}\right)=\int_{0}^{\infty}\frac{\cos\left(\frac{\pi x}{2}\right)}{\{-1+\exp{(2\pi\sqrt{x})}\}}dx=\frac{1}{4\pi},
\end{equation}
\begin{equation}\label{RCG15}
\Phi_3\left(\frac{2}{5}\right)=\int_{0}^{\infty}\frac{\cos\left(\frac{2\pi x}{5}\right)}{\{-1+\exp{(2\pi\sqrt{x})}\}}dx=\frac{8-3\sqrt{5}}{16}.
\end{equation}
From eqns(\ref{CR27}) and (\ref{CR29}), we obtain
\begin{equation}\label{SR19}
\Psi_{3}^{*}(n)=\int_{0}^{\infty}\frac{\sin(\pi nx)}{\{-1+\exp{(2\pi\sqrt{x})}\}}dx=\frac{1}{n}\sqrt{\left(\frac{2}{n}\right)}~\Phi_{3}\left(\frac{1}{n}\right)+\Phi_{3}(n)-\frac{1}{2\pi n}.
\end{equation}
Setting $n=1,2,\frac{1}{2}$ in the above eq. (\ref{SR19}),  using values of $\Phi_{3}(1),\Phi_{3}(2)$ and $\Phi_{3}\left(\frac{1}{2}\right)$ from the eqns (\ref{RCG10}),(\ref{RCG11}) and (\ref{RCG14}), after simplification we get the following three results:
\begin{equation}\label{SR9}
\Psi_{3}^{*}(1)=\int_{0}^{\infty}\frac{\sin(\pi x)}{\{-1+\exp{(2\pi\sqrt{x})}\}}dx=\frac{\pi\sqrt{2}-4}{8\pi},
\end{equation}
\begin{equation}\label{SR10}
\Psi_{3}^{*}(2)=\int_{0}^{\infty}\frac{\sin(2\pi x)}{\{-1+\exp{(2\pi\sqrt{x})}\}}dx=\frac{\pi-2}{16\pi},
\end{equation}
\begin{equation}\label{SR11}
\Psi_{3}^{*}(1/2)=\int_{0}^{\infty}\frac{\sin\left(\frac{\pi x}{2}\right)}{\{-1+\exp{(2\pi\sqrt{x})}\}}dx=\frac{\pi-3}{4\pi}.
\end{equation}
Binomial function is given by
\begin{equation}\label{FAR8}
(1-z)^{-a}={}_{1}F_{0}\left(\begin{array}{lll}a~;\\  \overline{~~~};\end{array} z \right)
=\sum_{n=0}^{\infty}\frac{(a)_{n}}{n!}z^{n},
\end{equation}
~~~~~where $|z|<1,~~a\in\mathbb{C}$.\\

The Meijer's $G$ function is defined by means of  the Mellin-Barnes type contour integral \cite[Sec.(1.5), eq. (1)]{s1}.
 When $k=1,2,...,n$ and $\ell=1,2,...,m$, and $\alpha_{k}-\beta_{\ell}\neq$ positive integer then
\begin{multline}\label{RCG32}
G_{{p} ,{q}}^{{m},{n}} \left( z~\bigg{|} \begin{array}{lll} \alpha_{1},...,\alpha_{n};\alpha_{n+1},...,\alpha_{p} \\\beta_{1},...,\beta_{m};\beta_{m+1},...,\beta{q} \end{array} \right)\\
=\frac{1}{2\pi i}\bigint_{-i\infty}^{+i\infty}
\frac{\displaystyle\prod_{j=1}^{m}\Gamma(\beta_{j}-\zeta)\prod_{j=1}^{n}\Gamma(1-\alpha_{j}+\zeta)}{\displaystyle\prod_{j=m+1}^{q}\Gamma(1-\beta_{j}+\zeta)\prod_{j=n+1}^{p}\Gamma(\alpha_{j}-\zeta)}(z)^{\zeta}d\zeta~,~~~~~~~~~~~~~~
\end{multline}
where $z\neq0$, and $m,n,p,q$ are non negative integers such that $ 1\leq m\leq q$~;~$0\leq n\leq p$.\\
Suppose:
\begin{equation}\label{RCG33}
\bigwedge=m+n-\frac{1}{2}(p+q),
\end{equation}
~and
\begin{equation}
  \omega=(\beta_{1}+...+\beta_{m}+\beta_{m+1}+...+\beta{q})-(\alpha_{1}+...+\alpha_{n}+\alpha_{n+1}+...+\alpha_{p}).
\end{equation}
(i) If $|arg(z)|<\bigwedge\pi$ and $\bigwedge > 0$, then the integral (\ref{RCG32}) is converges.\\
(ii) If $|arg(z)|=\bigwedge\pi$ and $\bigwedge\geq0$, then the integral (\ref{RCG32}) converges absolutely\\ when $p=q$ and $\Re(\omega)<-1$.\\
(iii) If $|arg(z)|=\bigwedge\pi$ and $\bigwedge\geq0$, then the integral (\ref{RCG32}) also converges absolutely\\ when $p\neq q$ and
\begin{equation}
  (q-p)\xi>1-\left(\frac{q-p}{2}\right)+\Re(\omega),
\end{equation}
where $\zeta=\xi+i\eta;\xi$ is so chosen that, for $\eta\rightarrow\pm\infty;\xi$ and $\eta$ are real.\\
Elementary property of G-function is given by \cite[p.46, eq. (6)]{s1}
\begin{multline}\label{PG1}
G_{{p} ,{q}}^{{m},{n}} \left( z~\bigg{|} \begin{array}{lll} \alpha_{1},...,\alpha_{n};\alpha_{n+1},...,\alpha_{p} \\\beta_{1},...,\beta_{m};\beta_{m+1},...,\beta{q} \end{array} \right)\\
=G_{{q} ,{p}}^{{n},{m}} \left( \frac{1}{z}~\bigg{|} \begin{array}{lll} 1-\beta_{1},...,1-\beta_{m};1-\beta_{m+1},...,1-\beta{q} \\1-\alpha_{1},...,1-\alpha_{n}; 1-\alpha_{n+1},...,1-\alpha_{p}  \end{array} \right),
\end{multline}
when $p>q$, then we can apply above property (\ref{PG1}) through out the paper.\\
The representation of sine function in terms of  G-Function \cite{E1, MA1, MA2, SG1}, is given by
\begin{equation}\label{RG38}
\sin(z)= \sqrt(\pi)~~G_{{0} ,{2}}^{{1},{0}} \left( \frac{z^{2}}{4}|\begin{array}{lll}\overline{~~~~~~~~~~~~~}\\ \frac{1}{2}~~;~~0\end{array} \right),
\end{equation}
\begin{equation}\label{RMG43}
=\sqrt{\pi}\frac{1}{2\pi i}\int_{-i\infty}^{+i\infty}\frac{\Gamma(\frac{1}{2}-\zeta)}{\Gamma(1+\zeta)}\left(\frac{z^{2}}{4}\right)^{\zeta}d\zeta.
\end{equation}
The above Mellin-Barnes type contour integral representation (\ref{RMG43}) of sine function is obtained by using the definition (\ref{RCG32}) in the equation (\ref{RG38}).\\
The representation of cosine function in terms of  G-Function \cite{E1, MA1, MA2, SG1} is given by
\begin{equation}\label{RCG38}
\cos(z)= \sqrt(\pi)~~G_{{0} ,{2}}^{{1},{0}} \left(\ \frac{z^{2}}{4}\bigg{|} \begin{array}{lll}\overline{~~~~~~~~~~~~~}\\ 0~~;~~ \frac{1}{2}\end{array} \right),
\end{equation}
\begin{equation}\label{RGC43}
=\sqrt{\pi}\frac{1}{2\pi i}\int_{-i\infty}^{+i\infty}\frac{\Gamma(-\zeta)}{\Gamma(\frac{1}{2}+\zeta)}\left(\frac{z^{2}}{4}\right)^{\zeta}d\zeta.
\end{equation}
The above Mellin-Barnes type contour integral representation (\ref{RGC43}) of cosine function  is obtained by using the definition (\ref{RCG32}) in the equation (\ref{RCG38}).\\
For every positive integer $m$ \cite[p.22, Eq.(26)]{s1}, we have
\begin{equation}\label{RCG50}
  \Gamma(mz)=(2\pi)^{\frac{(1-m)}{2}}m^{mz-\frac{1}{2}}\prod_{j=1}^{m}\Gamma\left(z+\frac{j-1}{m}\right),~~ mz\in\mathbb{C}\setminus \mathbb{Z}_{0}^{-}.
\end{equation}
The equation (\ref{RCG50}) is known as Gauss-Legendre multiplication theorem for Gamma function.\\

Laplace transform of $t^{z-1}$ \cite[p.12, eq.(33)]{E1} is given by
\begin{equation}\label{FAR23}
\mathcal{L}[t^{z-1};S]=\int_{0}^{\infty}~e^{-St}t^{z-1}dt=\frac{\Gamma(z)}{S^{z}},
\end{equation}
~~~~where $ \Re(S)>0 ~and~ 0<\Re(z)<\infty $.\\
The Laplace transforms of sine and cosine functions associated with Meijer's $G$-function are given by
\begin{equation}\label{G1}
\mathcal{L}[\sin(\beta x^{2}]=\int_{0}^{\infty}~e^{-\alpha x} \sin(\beta x^{2})dx=\frac{1}{\alpha\pi\sqrt{2}}
~G_{{3} ,{1}}^{{1},{3}} \left(\ \frac{64\beta^{2}}{\alpha^{4}}\bigg{|} \begin{array}{lll}\frac{1}{4},\frac{1}{2},\frac{3}{4}\\ \frac{1}{2} \end{array} \right),
\end{equation}
\begin{equation}\label{G2}
\mathcal{L}[\cos(\beta x^{2})]=\int_{0}^{\infty}~e^{-\alpha x} \cos(\beta x^{2})dx=\frac{1}{\alpha\pi\sqrt{2}}
~G_{{3} ,{1}}^{{1},{3}} \left(\ \frac{64\beta^{2}}{\alpha^{4}}\bigg{|} \begin{array}{lll}\frac{1}{4},\frac{3}{4},0\\0\end{array} \right),
\end{equation}
\begin{equation}\label{G3}
\mathcal{L}[x\sin(\beta x^{2})]=\int_{0}^{\infty}~e^{-\alpha x} x~\sin(\beta x^{2})dx=\frac{2\sqrt{2}}{\alpha^{2}\pi}
~G_{{3} ,{1}}^{{1},{3}} \left(\ \frac{64\beta^{2}}{\alpha^{4}}\bigg{|} \begin{array}{lll}-\frac{1}{4},\frac{1}{4},\frac{1}{2}\\ \frac{1}{2} \end{array} \right),
\end{equation}
\begin{equation}\label{G4}
\mathcal{L}[x\cos(\beta x^{2})]=\int_{0}^{\infty}~e^{-\alpha x}x~ \cos(\beta x^{2})dx=\frac{2\sqrt{2}}{\alpha^{2}\pi}
~G_{{3} ,{1}}^{{1},{3}} \left(\ \frac{64\beta^{2}}{\alpha^{4}}\bigg{|} \begin{array}{lll}-\frac{1}{4},\frac{1}{4},0\\0 \end{array} \right),
\end{equation}
~~~~~~~~where $\Re(\alpha)>0$.\\
 \textbf{Proof}: The formulas (\ref{G1})-(\ref{G4})  are obtained by using the  Mellin-Barnes type contour integrals (\ref{RMG43}) and (\ref{RGC43}) of  sine and cosine functions.  Then change the order of integration in double integrals, take Laplace transform with the help of eq.(\ref{FAR23})and apply  Gauss-Legendre multiplication theorem (\ref{RCG50}) for the factors of $\Gamma (4\zeta+1), ~\Gamma (4\zeta+2)$  and finally use the definition of Meijer's G-function (\ref{RCG32}), we get above results (\ref{G1})-(\ref{G4}); see also in \cite{E2}. \\
 \\
$~~~~~~~~$In this paper our integrals (\ref{G1})-(\ref{G4}) play an important role in our investigation. Using Laplace transforms formulas (\ref{G1})-(\ref{G4}), we solved all Ramanujan's integrals involved in three theorems. Further, we obtain six infinite summation formulas involving Meijer's G-function. Numeric values of thirteen infinite series are also discussed.
\section{ Analytical solutions of Ramanujan's integrals}
\texttt{Each of the following  Ramanujan's integrals involving infinite series of \\Meijer's $G$-function holds true}:\\
\texttt{\textbf{I}. The first Ramanujan integrals holds true}
\begin{equation}\label{RG11}
\Psi_{1}(b)=\int_{0}^{\infty}\frac{\sin(b\pi x^{2})}{\cosh(\pi x)}dx
=\frac{\sqrt{2}}{\pi^{2} }\sum_{r=0}^{\infty}\bigg[\frac{(-1)^{r}}{(1+2r)}
~G_{{3} ,{1}}^{{1},{3}} \left(\ \frac{64b^{2}}{\pi^{2}(1+2r)^{4}}\bigg{|} \begin{array}{lll}\frac{1}{4},\frac{1}{2},\frac{3}{4}\\ \frac{1}{2}\end{array} \right)\bigg],
\end{equation}
\begin{equation}\label{RG12}
\Phi_{1}(b)=\int_{0}^{\infty}\frac{\cos(b\pi x^{2})}{\cosh(\pi x)}dx
=\frac{\sqrt{2}}{\pi^{2} }\sum_{r=0}^{\infty}\bigg[\frac{(-1)^{r}}{(1+2r)}
~G_{{3} ,{1}}^{{1},{3}} \left(\ \frac{64b^{2}}{\pi^{2}(1+2r)^{4}}\bigg{|} \begin{array}{lll}\frac{1}{4},\frac{3}{4},0\\ 0\end{array} \right)\bigg].
\end{equation}
\texttt{\textbf{II}. The second Ramanujan integrals holds true}
\begin{multline}\label{RG13}
\Psi_{2}(b)=\int_{0}^{\infty}\frac{\sin(b\pi x^{2})}{\left\{1+2\cosh\left(\frac{2\pi x}{\sqrt{3}}\right)\right\}}dx\\
=\frac{\sqrt{3}}{\pi^{2}2\sqrt{2}}\sum_{p,q=0}^{\infty}\bigg[\frac{(1)_{p+q}(1)_{2q+p}(-1)^{p+q}}{(2)_{2q+p}~p!~q!}
~G_{{3} ,{1}}^{{1},{3}} \left(\ \frac{36b^{2}}{\pi^{2}}\left\{\frac{(1)_{2q+p}}{(2)_{2q+p}}\right\}^{4}\bigg{|} \begin{array}{lll}\frac{1}{4},\frac{1}{2},\frac{3}{4}\\ \frac{1}{2}\end{array} \right)\bigg],
\end{multline}
\begin{multline}\label{RG14}
\Phi_{2}(b)=\int_{0}^{\infty}\frac{\cos(b\pi x^{2})}{\left\{1+2\cosh\left(\frac{2\pi x}{\sqrt{3}}\right)\right\}}dx\\
=\frac{\sqrt{3}}{\pi^{2}2\sqrt{2}}\sum_{p,q=0}^{\infty}\bigg[\frac{(1)_{p+q}(1)_{2q+p}(-1)^{p+q}}{(2)_{2q+p}~p!~q!}
~G_{{3} ,{1}}^{{1},{3}} \left(\ \frac{36b^{2}}{\pi^{2}}\left\{\frac{(1)_{2q+p}}{(2)_{2q+p}}\right\}^{4}\bigg{|} \begin{array}{lll}\frac{1}{4},\frac{3}{4},0\\ 0\end{array} \right)\bigg].
\end{multline}
\texttt{\textbf{III}. The third Ramanujan integrals holds true}
\begin{equation}\label{RG15}
\Psi_{3}^{*}(b)=\int_{0}^{\infty}\frac{\sin(b\pi x)}{\{-1+\exp{(2\pi\sqrt{x})}\}}dx
=\frac{\sqrt{2}}{\pi^{3}}\sum_{r=0}^{\infty}\bigg[\left\{\frac{(1)_{r}}{(2)_{r}}\right\}^{2}
~G_{{3} ,{1}}^{{1},{3}} \left(\ \frac{4b^{2}}{\pi^{2}}\left\{\frac{(1)_{r}}{(2)_{r}}\right\}^{4}\bigg{|} \begin{array}{lll}-\frac{1}{4},\frac{1}{4},\frac{1}{2}\\ \frac{1}{2}\end{array} \right)\bigg],
\end{equation}
\begin{equation}\label{RG16}
\Phi_{3}(b)=\int_{0}^{\infty}\frac{\cos(b\pi x)}{\{-1+\exp{(2\pi\sqrt{x})}\}}dx
=\frac{\sqrt{2}}{\pi^{3}}\sum_{r=0}^{\infty}\bigg[\left\{\frac{(1)_{r}}{(2)_{r}}\right\}^{2}
~G_{{3} ,{1}}^{{1},{3}} \left(\ \frac{4b^{2}}{\pi^{2}}\left\{\frac{(1)_{r}}{(2)_{r}}\right\}^{4}\bigg{|} \begin{array}{lll}-\frac{1}{4},\frac{1}{4},0\\ 0\end{array} \right)\bigg].
\end{equation}
\textbf{Proof}:  Suppose left hand side of the eq.(\ref{RG11}) is denoted by $\Psi_{1}(b)$ upon using the well known results of hyperbolic functions and  Binomial function (\ref{FAR8}), we obtain
\begin{eqnarray}\label{AF06}
\Psi_{1}(b)=2\int_{0}^{\infty}e^{-\pi x}(1+e^{-2\pi x})^{-1}\sin(b\pi x^{2})dx,\\
=2\int_{0}^{\infty}e^{-\pi x}{}_{1}F_{0}\left(\begin{array}{lll}1~;\\  \overline{~~~};\end{array} -e^{-2\pi x} \right)\sin(b\pi x^{2})dx.
 \end{eqnarray}
 Change the order of integration and summation in above equation, which yield
 \begin{equation}\label{A006}
 \Psi_{1}(b)=2\sum_{r=0}^{\infty}(-1)^{r}\int_{0}^{\infty}e^{-(\pi+2\pi r)x}\sin(b\pi x^{2})dx.
 \end{equation}
 Use Laplace formula (\ref{G1}) in the eq.(\ref{A006}). Then we get the right hand side stated in (\ref{RG11}). Similarly, proof of the  integral (\ref{RG12})  by using the Laplace formula (\ref{G2}) is much akin as result (\ref{RG11}), which we have already discussed in a detailed manner.

 Again suppose left hand side of the eq.(\ref{RG13}) is denoted by $\Psi_{2}(b)$, we have
 \begin{eqnarray}\label{RR01}
\Psi_{2}(b)=\int_{0}^{\infty}\frac{\sin(b\pi x^{2})}{(1+e^{2\pi x/\sqrt{3}})+e^{-2\pi x/\sqrt{3}}}dx,
 \end{eqnarray}
 \begin{eqnarray}\label{RR03}
=\int_{0}^{\infty}\frac{1}{(1+e^{2\pi x/\sqrt{3}})}\left\{1+\frac{e^{-2\pi x/\sqrt{3}}}{1+e^{2\pi x/\sqrt{3}}}\right\}^{-1}\sin(b\pi x^{2})dx.
 \end{eqnarray}
 Change the order of integration and summation in above equation, which yield
\begin{eqnarray}\label{RR04}
\Psi_{2}(b)=\sum_{q=0}^{\infty}(-1)^{q}\int_{0}^{\infty}\frac{e^{-2\pi qx/\sqrt{3}}}{(1+e^{2\pi x/\sqrt{3}})^{q+1}}\sin(b\pi x^{2})dx,
 \end{eqnarray}
\begin{eqnarray}\label{RR05}
=\sum_{q=0}^{\infty}(-1)^{q}\int_{0}^{\infty} {}_{1}F_{0}\left(\begin{array}{lll}q+1~;\\  \overline{~~~};\end{array} -e^{-2\pi x/\sqrt{3}} \right)e^{-(4\pi q+2\pi)x/\sqrt{3}}\sin(b\pi x^{2})dx,
 \end{eqnarray}
\begin{eqnarray}\label{RR06}
=\sum_{p=0}^{\infty}\sum_{q=0}^{\infty}\frac{(-1)^{p+q}(1)_{p+q}}{p! q!}\int_{0}^{\infty} e^{-\frac{(4\pi q+2\pi p+2\pi)x}{\sqrt{3}}}\sin(b\pi x^{2})dx.
 \end{eqnarray}
 Use Laplace formula (\ref{G1}) in the eq.(\ref{RR06}). Then we get the right hand side stated in (\ref{RG13}). Similarly, proof of the  integral (\ref{RG14}) by using the Laplace formulas (\ref{G2}) is much akin as result (\ref{RG13}), which we have already discussed in a detailed manner.

Suppose left hand side of the eq.(\ref{RG15}) is denoted by $\Psi_{3}^{*}(b)$, we have
\begin{eqnarray}\label{RR07}
\Psi_{3}^{*}(b)=\int_{0}^{\infty}e^{-2\pi \sqrt{x}}(1-e^{-2\pi \sqrt{x}})^{-1}\sin(b\pi x)dx,
\end{eqnarray}
\begin{eqnarray}
=\int_{0}^{\infty}e^{-2\pi \sqrt{x}}{}_{1}F_{0}\left(\begin{array}{lll}1~;\\  \overline{~~~};\end{array} e^{-2\pi \sqrt{x}} \right)\sin(b\pi x)dx.
 \end{eqnarray}
 Change the order of integration and summation in above equation, which yield
 \begin{equation}\label{RR07}
 \Psi_{3}^{*}(b)=\sum_{r=0}^{\infty}\int_{0}^{\infty}e^{-(2\pi+2\pi r)\sqrt{x}}\sin(b\pi x)dx.
 \end{equation}
 Setting $\sqrt{x}=t$, when $x = 0,~~t = 0$, when $x\rightarrow\infty,~~t\rightarrow\infty$, in the L.H.S. of eq.(\ref{RR07}), we obtain
 \begin{equation}\label{RR08}
 \Psi_{3}^{*}(b)=2\sum_{r=0}^{\infty}\int_{0}^{\infty}e^{-(2\pi+2\pi r)t}t\sin(b\pi t^{2})dt.
 \end{equation}
 Use Laplace formula (\ref{G3}) in the eq.(\ref{RR08}). Then we get the right hand side stated in eq.(\ref{RG15}). Similarly, proof of the  integral (\ref{RG16}) by using the Laplace formula (\ref{G4}) is much akin as result (\ref{RG15}), which we have already discussed in a detailed manner.
\section{Application of Ramanujan's integrals in three theorems}
The following infinite summation formulas involving Meijer's G-function hold true.

\begin{multline}\label{RG31}
\sum_{r=0}^{\infty}\bigg[\frac{(-1)^{r}}{(1+2r)}
~G_{{3} ,{1}}^{{1},{3}} \left(\ \frac{64n^2}{\pi^{2}(1+2r)^{4}}\bigg{|} \begin{array}{lll}\frac{1}{4},\frac{3}{4},0\\0\end{array} \right)\bigg]\\
=\sum_{r=0}^{\infty}\frac{(-1)^{r}}{(1+2r)}
~\bigg[\sqrt{\frac{2}{n}}~~G_{{3} ,{1}}^{{1},{3}} \left(\ \frac{64}{(n\pi)^{2}(1+2r)^{4}}\bigg{|} \begin{array}{lll}\frac{1}{4},\frac{1}{2},\frac{3}{4}\\ \frac{1}{2}\end{array} \right)\\
+G_{{3} ,{1}}^{{1},{3}} \left(\ \frac{64n^2}{\pi^{2}(1+2r)^{4}}\bigg{|} \begin{array}{lll}\frac{1}{4},\frac{1}{2},\frac{3}{4}\\ \frac{1}{2}\end{array} \right)\bigg],
\end{multline}
\begin{multline}\label{RG32}
\sum_{r=0}^{\infty}\bigg[\frac{(-1)^{r}}{(1+2r)}
~G_{{3} ,{1}}^{{1},{3}} \left(\ \frac{64n^2}{\pi^{2}(1+2r)^{4}}\bigg{|} \begin{array}{lll}\frac{1}{4},\frac{1}{2},\frac{3}{4}\\ \frac{1}{2}\end{array} \right)\bigg]\\
=\sum_{r=0}^{\infty}\frac{(-1)^{r}}{(1+2r)}
~\bigg[\sqrt{\frac{2}{n}}~~G_{{3} ,{1}}^{{1},{3}} \left(\ \frac{64}{(n\pi)^{2}(1+2r)^{4}}\bigg{|} \begin{array}{lll}\frac{1}{4},\frac{3}{4},0\\ 0\end{array} \right)\\
-G_{{3} ,{1}}^{{1},{3}} \left(\ \frac{64n^2}{\pi^{2}(1+2r)^{4}}\bigg{|} \begin{array}{lll}\frac{1}{4},\frac{3}{4},0\\ 0\end{array} \right)\bigg],
\end{multline}
\begin{multline}\label{RG33}
\sum_{p,q=0}^{\infty}\bigg[\frac{(1)_{p+q}(1)_{2q+p}(-1)^{p+q}}{(2)_{2q+p}~p! q!}
~G_{{3} ,{1}}^{{1},{3}} \left(\ \frac{36n^2}{\pi^{2}}\left\{\frac{(1)_{2q+p}}{(2)_{2q+p}}\right\}^{4}\bigg{|} \begin{array}{lll}\frac{1}{4},\frac{3}{4},0\\0\end{array} \right)\bigg]\\
=\sum_{p,q=0}^{\infty}\frac{(1)_{p+q}(1)_{2q+p}(-1)^{p+q}}{(2)_{2q+p}~p! q!}
~\bigg[\sqrt{\frac{2}{n}}~~G_{{3} ,{1}}^{{1},{3}} \left(\ \frac{36}{(n\pi)^{2}}\left\{\frac{(1)_{2q+p}}{(2)_{2q+p}}\right\}^{4}\bigg{|} \begin{array}{lll}\frac{1}{4},\frac{1}{2},\frac{3}{4}\\ \frac{1}{2}\end{array} \right)\\
+G_{{3} ,{1}}^{{1},{3}} \left(\ \frac{36n^2}{\pi^{2}}\left\{\frac{(1)_{2q+p}}{(2)_{2q+p}}\right\}^{4}\bigg{|} \begin{array}{lll}\frac{1}{4},\frac{1}{2},\frac{3}{4}\\ \frac{1}{2}\end{array} \right)\bigg],
\end{multline}
\begin{multline}\label{RG34}
\sum_{p,q=0}^{\infty}\bigg[\frac{(1)_{p+q}(1)_{2q+p}(-1)^{p+q}}{(2)_{2q+p}~p! q!}
~G_{{3} ,{1}}^{{1},{3}} \left(\ \frac{36n^2}{\pi^{2}}\left\{\frac{(1)_{2q+p}}{(2)_{2q+p}}\right\}^{4}\bigg{|} \begin{array}{lll}\frac{1}{4},\frac{1}{2},\frac{3}{4}\\ \frac{1}{2}\end{array} \right)\bigg]\\
=\sum_{p,q=0}^{\infty}\frac{(1)_{p+q}(1)_{2q+p}(-1)^{p+q}}{(2)_{2q+p}~p! q!}
~\bigg[\sqrt{\frac{2}{n}}~~G_{{3} ,{1}}^{{1},{3}} \left(\ \frac{36}{(n\pi)^{2}}\left\{\frac{(1)_{2q+p}}{(2)_{2q+p}}\right\}^{4}\bigg{|} \begin{array}{lll}\frac{1}{4},\frac{3}{4},0\\ 0\end{array} \right)\\
-G_{{3} ,{1}}^{{1},{3}} \left(\ \frac{36n^2}{\pi^{2}}\left\{\frac{(1)_{2q+p}}{(2)_{2q+p}}\right\}^{4}\bigg{|} \begin{array}{lll}\frac{1}{4},\frac{3}{4},0\\ 0\end{array} \right)\bigg],
\end{multline}
\begin{multline}\label{RG35}
\sum_{r=0}^{\infty}\bigg[\left\{\frac{(1)_{r}}{(2)_{r}}\right\}^{2}
~G_{{3} ,{1}}^{{1},{3}} \left(\ \frac{4n^2}{\pi^{2}}\left\{\frac{(1)_{r}}{(2)_{r}}\right\}^{4}\bigg{|} \begin{array}{lll}-\frac{1}{4},\frac{1}{4},\frac{1}{2}\\ \frac{1}{2}\end{array} \right)\bigg]\\
=-\frac{\pi^{2}}{2n\sqrt{2}}+\sum_{r=0}^{\infty}\left\{\frac{(1)_{r}}{(2)_{r}}\right\}^{2}
~\bigg[\frac{\sqrt{2}}{n\sqrt{n}}~~G_{{3} ,{1}}^{{1},{3}} \left(\ \frac{4}{(n\pi)^{2}}\left\{\frac{(1)_{r}}{(2)_{r}}\right\}^{4}\bigg{|} \begin{array}{lll}-\frac{1}{4},\frac{1}{4},0\\ 0\end{array} \right)\\
+G_{{3} ,{1}}^{{1},{3}} \left(\ \frac{4n^2}{\pi^{2}}\left\{\frac{(1)_{r}}{(2)_{r}}\right\}^{4}\bigg{|} \begin{array}{lll}-\frac{1}{4},\frac{1}{4},0\\ 0\end{array} \right)\bigg],
\end{multline}
\begin{multline}\label{RG36}
\sum_{r=0}^{\infty}\bigg[\left\{\frac{(1)_{r}}{(2)_{r}}\right\}^{2}
~G_{{3} ,{1}}^{{1},{3}} \left(\ \frac{4n^2}{\pi^{2}}\left\{\frac{(1)_{r}}{(2)_{r}}\right\}^{4}\bigg{|} \begin{array}{lll}-\frac{1}{4},\frac{1}{4},0\\ 0\end{array} \right)\bigg]\\
=\frac{\pi^{2}}{2}\left(\frac{1}{\sqrt{n}}-\frac{1}{n\sqrt{2}}\right)+\sum_{r=0}^{\infty}\left\{\frac{(1)_{r}}{(2)_{r}}\right\}^{2}
~\bigg[\frac{\sqrt{2}}{n\sqrt{n}}~~G_{{3} ,{1}}^{{1},{3}} \left(\ \frac{4}{(n\pi)^{2}}\left\{\frac{(1)_{r}}{(2)_{r}}\right\}^{4}\bigg{|} \begin{array}{lll}-\frac{1}{4},\frac{1}{4},\frac{1}{2}\\ \frac{1}{2}\end{array} \right)\\
-G_{{3} ,{1}}^{{1},{3}} \left(\ \frac{4n^2}{\pi^{2}}\left\{\frac{(1)_{r}}{(2)_{r}}\right\}^{4}\bigg{|} \begin{array}{lll}-\frac{1}{4},\frac{1}{4},\frac{1}{2}\\ \frac{1}{2} \end{array} \right)\bigg].
\end{multline}

\textbf{Proof}: Put the values of $\Phi_1(n),\Phi_1(1/n);\Psi_1(n),\Psi_1(1/n);\Phi_2(n),\Phi_2(1/n);\Psi_2(n),\Psi_2(1/n);\\
\Phi_3(n),\Phi_3(1/n);\Psi_3(n),\Psi_3(1/n); \Psi_{3}^{*}(n),\Psi_{3}^{*}(1/n)$ with the help of Ramanujan's integrals (\ref{RG11})-(\ref{RG16}) in the equations (\ref{CR16}),(\ref{CR17}),(\ref{CR24}),(\ref{CR25}),(\ref{CR28}),(\ref{CR29}), after simplification. Then we get six infinite summation formulas (\ref{RG31})-(\ref{RG36}).

\section{Numerical values of some infinite series containing Meijer's $G$-function}
The numerical values of thirteen infinite series associated  with Meijer's $G$-function hold true

\begin{equation}\label{SF1}
\sum_{r=0}^{\infty}\bigg[\frac{(-1)^{r}}{(1+2r)}
~G_{{3} ,{1}}^{{1},{3}} \left(\ \frac{64}{\pi^{2}(1+2r)^{4}}\bigg{|} \begin{array}{lll}\frac{1}{4},\frac{1}{2},\frac{3}{4}\\ \frac{1}{2} \end{array} \right)\bigg]=\frac{\pi^{2}}{4}(-1+\sqrt{2}),
\end{equation}
\begin{equation}\label{SF2}
\sum_{r=0}^{\infty}\bigg[\frac{(-1)^{r}}{(1+2r)}
~G_{{3} ,{1}}^{{1},{3}} \left(\ \frac{64}{\pi^{2}(1+2r)^{4}}\bigg{|} \begin{array}{lll}\frac{1}{4},\frac{3}{4},0\\ 0 \end{array} \right)\bigg]=\frac{\pi^{2}}{4},
\end{equation}
\begin{multline}\label{SF3}
\sum_{p,q=0}^{\infty}\bigg[\frac{(1)_{p+q}(1)_{2q+p}(-1)^{p+q}}{(2)_{2q+p}~p! q!}
~G_{{3} ,{1}}^{{1},{3}} \left(\ \frac{36}{\pi^{2}}\left\{\frac{(1)_{2q+p}}{(2)_{2q+p}}\right\}^{4}\bigg{|} \begin{array}{lll}\frac{1}{4},\frac{1}{2},\frac{3}{4}\\ \frac{1}{2}\end{array} \right)\bigg]\\=\frac{\pi^{2}}{2\sqrt{3}}(-\sqrt{6}+\sqrt{3}+4),
\end{multline}
\begin{multline}\label{SF4}
\sum_{p,q=0}^{\infty}\bigg[\frac{(1)_{p+q}(1)_{2q+p}(-1)^{p+q}}{(2)_{2q+p}~p! q!}
~G_{{3} ,{1}}^{{1},{3}} \left(\ \frac{36}{\pi^{2}}\left\{\frac{(1)_{2q+p}}{(2)_{2q+p}}\right\}^{4}\bigg{|} \begin{array}{lll}\frac{1}{4},\frac{3}{4},0\\ 0\end{array} \right)\bigg]\\=\frac{\pi^{3}}{2\sqrt{3}}(\sqrt{2}-\sqrt{3}+1),
\end{multline}
\begin{equation}\label{NG1}
\sum_{r=0}^{\infty}\bigg[\left\{\frac{(1)_{r}}{(2)_{r}}\right\}^{2}
~G_{{3} ,{1}}^{{1},{3}} \left(\ \frac{4}{\pi^{2}}\left\{\frac{(1)_{r}}{(2)_{r}}\right\}^{4}\bigg{|} \begin{array}{lll}-\frac{1}{4},\frac{1}{4},\frac{1}{2}\\ \frac{1}{2}\end{array} \right)\bigg]=\frac{\pi^{2}}{8}(\pi-2\sqrt{2}),
\end{equation}
\begin{equation}\label{NG2}
\sum_{r=0}^{\infty}\bigg[\left\{\frac{(1)_{r}}{(2)_{r}}\right\}^{2}
~G_{{3} ,{1}}^{{1},{3}} \left(\ \frac{16}{\pi^{2}}\left\{\frac{(1)_{r}}{(2)_{r}}\right\}^{4}\bigg{|} \begin{array}{lll}-\frac{1}{4},\frac{1}{4},\frac{1}{2}\\ \frac{1}{2}\end{array} \right)\bigg]=\frac{\pi^{2}\sqrt{2}}{32}(\pi-2),
\end{equation}
\begin{equation}\label{NG3}
\sum_{r=0}^{\infty}\bigg[\left\{\frac{(1)_{r}}{(2)_{r}}\right\}^{2}
~G_{{3} ,{1}}^{{1},{3}} \left(\ \frac{1}{\pi^{2}}\left\{\frac{(1)_{r}}{(2)_{r}}\right\}^{4}\bigg{|} \begin{array}{lll}-\frac{1}{4},\frac{1}{4},\frac{1}{2}\\ \frac{1}{2}\end{array} \right)\bigg]=\frac{\pi^{2}\sqrt{2}}{8}(\pi-3),
\end{equation}
\begin{equation}\label{NG4}
\sum_{r=0}^{\infty}\bigg[\left\{\frac{(1)_{r}}{(2)_{r}}\right\}^{2}
~G_{{3} ,{1}}^{{1},{3}} \left(\ \frac{4}{\pi^{2}}\left\{\frac{(1)_{r}}{(2)_{r}}\right\}^{4}\bigg{|} \begin{array}{lll}-\frac{1}{4},\frac{1}{4},0\\ 0\end{array} \right)\bigg]=\frac{\pi^{3}}{16}(2\sqrt{2}-2),
\end{equation}
\begin{equation}\label{RG40}
\sum_{r=0}^{\infty}\bigg[\left\{\frac{(1)_{r}}{(2)_{r}}\right\}^{2}
~G_{{3} ,{1}}^{{1},{3}} \left(\ \frac{16}{\pi^{2}}\left\{\frac{(1)_{r}}{(2)_{r}}\right\}^{4}\bigg{|} \begin{array}{lll}-\frac{1}{4},\frac{1}{4},0\\ 0\end{array} \right)\bigg]=\frac{\pi^{3}\sqrt{3}}{48},
\end{equation}
\begin{equation}\label{RG41}
\sum_{r=0}^{\infty}\bigg[\left\{\frac{(1)_{r}}{(2)_{r}}\right\}^{2}
~G_{{3} ,{1}}^{{1},{3}} \left(\ \frac{1}{\pi^{2}}\left\{\frac{(1)_{r}}{(2)_{r}}\right\}^{4}\bigg{|} \begin{array}{lll}-\frac{1}{4},\frac{1}{4},0\\ 0\end{array} \right)\bigg]=\frac{\pi^{2}\sqrt{2}}{8},
\end{equation}
\begin{equation}\label{RG42}
\sum_{r=0}^{\infty}\bigg[\left\{\frac{(1)_{r}}{(2)_{r}}\right\}^{2}
~G_{{3} ,{1}}^{{1},{3}} \left(\ \frac{64}{\pi^{2}}\left\{\frac{(1)_{r}}{(2)_{r}}\right\}^{4}\bigg{|} \begin{array}{lll}-\frac{1}{4},\frac{1}{4},0\\ 0\end{array} \right)\bigg]=\frac{\pi^{3}}{64}(3\sqrt{2}-2),
\end{equation}
\begin{equation}\label{RG43}
\sum_{r=0}^{\infty}\bigg[\left\{\frac{(1)_{r}}{(2)_{r}}\right\}^{2}
~G_{{3} ,{1}}^{{1},{3}} \left(\ \frac{144}{\pi^{2}}\left\{\frac{(1)_{r}}{(2)_{r}}\right\}^{4}\bigg{|} \begin{array}{lll}-\frac{1}{4},\frac{1}{4},0\\ 0\end{array} \right)\bigg]=\frac{\pi^{3}}{288}(13\sqrt{2}-4\sqrt{6}),
\end{equation}
\begin{equation}\label{RG45}
\sum_{r=0}^{\infty}\bigg[\left\{\frac{(1)_{r}}{(2)_{r}}\right\}^{2}
~G_{{3} ,{1}}^{{1},{3}} \left(\ \frac{16}{25\pi^{2}}\left\{\frac{(1)_{r}}{(2)_{r}}\right\}^{4}\bigg{|} \begin{array}{lll}-\frac{1}{4},\frac{1}{4},0\\ 0\end{array} \right)\bigg]=\frac{\pi^{3}}{38}(8\sqrt{2}-3\sqrt{10}).
\end{equation}
\textbf{Proof}: Taking a particular value $b=1$ in the Ramanujan integrals (\ref{RG11})-(\ref{RG14}) and comparing with the corresponding values of the integrals (\ref{N1})-(\ref{N4}), then we get some numeric values of infinite series (\ref{SF1})-(\ref{SF4}). Again, taking some particular values of $b=1,2,1/2$ in the Ramanujan integral (\ref{RG15}) and comparing with values of the integrals (\ref{SR9})-(\ref{SR11}), then we get (\ref{NG1})-(\ref{NG3}).
 Similarly, taking some particular values of $b$ that is $b=1,2,4,6,1/2,2/5$  in the  Ramanujan integral (\ref{RG16}) and comparing with the values of integrals (\ref{RCG10})-(\ref{RCG15}), which gives some numeric values of infinite series (\ref{NG4})-(\ref{RG45}).
\section*{Concluding remarks}
  We have given the analytical solutions of  Srinivasa Ramanujan integrals in terms of Meijer's $G$-function. Further, we  have described some infinite summation formulas  connected with Meijer's G-function and  numeric values of some infinite series are also discussed. The numerical values of involved Ramanujan integrals are mentioned in our communicated paper \cite{Q13}, that supports to the corresponding our results of this paper. We conclude our present investigations by observing that several other similar integrals (\ref{RG11})-(\ref{RG16}) and consequences of infinite summation formulas of section-3 and its applications in numeric values of infinite series in section-4, can also be deduced in an analogous manner.

\end{document}